\documentclass{amsart}
\usepackage{latexsym,amssymb}
\usepackage{epsfig}
\usepackage{amsfonts}
\usepackage{amscd}
\usepackage{amsmath}
\theoremstyle{plain} \numberwithin{equation}{section}
\newtheorem{thm}{Theorem}[section]

\newtheorem{conj}{Conjecture}
\newtheorem{lem}[thm]{Lemma}

\theoremstyle{definition}
\newtheorem{de}[thm]{Definition}

\newtheorem{rem}[thm]{Remark}
\newtheorem{Example}[thm]{Example}

\newcommand{\F}{\mathbb{F}}
\newcommand{\R}{\mathbb{R}}
\newcommand{\C}{\mathbb{C}}

\newcommand{\Z}{\mathbb{Z}}

\newcommand{\cA}{\mathcal{A}}

\newcommand{\cD}{\mathcal{D}}
\newcommand{\cF}{\mathcal{F}}
\newcommand{\cG}{\mathcal{G}}
\newcommand{\cB}{\mathcal{B}}

\newcommand{\cM}{\mathcal{M}}
\newcommand{\cN}{\mathcal{N}}

\newcommand{\cL}{\mathcal{L}}
\newcommand{\cT}{\mathcal{T}}

\newcommand{\GL}{\mathrm{GL}}

\newcommand{\SO}{\mathrm{SO}}

\newcommand{\id}{\mathrm{id}}

\begin{document}

\begin{abstract}
Most of the examples of wavelet sets are for dilation sets
which are groups. We find a necessary and sufficient
condition under which subspace wavelet sets exist for
dilation sets of the form $\mathcal{A}\mathcal{B}$, which
are not necessarily groups. We explain the construction
by a few examples.
\end{abstract}

\author{Mihaela Dobrescu}
\email{dobrescu@math.lsu.edu}
\author{Gestur \'Olafsson}
\email{olafsson@math.lsu.edu}
\address{Both authors: Department of Mathematics\\
Louisiana State University\\
Baton Rouge, LA 70803, USA}

\title{Wavelet sets without groups}

\thanks{M. Dobrescu was partially
supported by DMS-0139783.
The research of G. \'Olafsson was supported by NSF grants DMS-0139783 and DMS-0402068.}
\maketitle
\section*{Introduction}
\noindent
Wavelets and frames have become a widely studied tool in mathematics and applied science.
One of the obvious questions is the construction of
wavelets with given properties. The classical
wavelet system on the line is given by a function  $\psi \in L^2(\R )$, such
that the dyadic dilates and translates of $\psi$ form an orthonormal
basis for $L^2(\R )$. Thus,
$\{\psi_{j,n}\}_{j,n\in\Z}$ with $\psi_{j,n}(t)=2^{j/2}\psi (2^jt+n)$, is
an orthonormal basis for $L^2(\R )$. There are several
obvious generalizations: One can replace $2$ by any integer $N$;
one can allow several wavelet functions $\psi^1,\ldots ,\psi^L$;
and one can consider an orthonormal basis for a closed subspace of
$L^2(\R )$. There have also been
several publications of wavelets in higher dimensions,
cf \cite{A03,ACDL98,BCMO95,BT96,FO02,F96,F98,HW89,LWWW2002,O02,olafsson_speegle:03}
to name few. One of the differences in higher dimensions is that
we now have many more choices in the sets of dilations and translations.
So, to fix the notation,
let $\mathcal {D}\subseteq GL(n,\mathbb{R})$ and
$\mathcal{T}\subseteq\mathbb{R}^n $ be countable sets.
A $(\mathcal D,\mathcal T)$-\textit{wavelet} is a square integrable
function $\psi$ with  the property that the set of functions
\begin{equation}\label{eq:wavelet}
\{|\det d|^{\frac{1}{2}}\varphi(dx+t)\mid d\in\mathcal{D},t\in\mathcal{T}
\}
\end{equation}
forms an orthonormal basis for $L^{2}(\mathbb{R}^{n})$.
The set $\mathcal {D}$ is  called the \textit{dilation set} and the set
$\mathcal{T}$ is called the \textit{translation set}.
If we replace  $L^{2}(\mathbb{R}^{n})$ in the above definition by
$$L^{2}_{M}(\mathbb{R}^{n})=\{f\in L^{2}(\mathbb{R}^{n})
\mid\mathrm{supp} (\mathcal{F}(f))\subseteq M\}$$
for some measurable subset $M\subseteq \mathbb{R}^n$, $|M|>0 $, we get a
$(\mathcal D,\mathcal T)$-{\textit subspace wavelet}.
Here $\cF$ stands
for the Fourier transform
$$\cF (f)(\lambda) =\int_{\R^n}f(t)e^{-2\pi i \langle{t , \lambda}
\rangle}\, dt\, .$$
We will often write $\hat{f}$ for the Fourier transform of $f$.

The most natural starting point is to consider groups of
dilations and full rank lattices as translation sets. The
simplest examples would then be groups generated by
one element $\cD=\{a^k\mid k\in \Z\}$,
see \cite{S} and
the reference therein.  In
\cite{O02,olafsson_speegle:03} more general sets of dilations were
considered, and in general those dilations do not form a group.
Even more general constructions can be found in \cite{A03}.

In this article, we consider a special class of wavelets corresponding to
\textit{wavelet sets}. Those are wavelet functions $\psi$ such that
$\cF (\psi)=\chi_{\Omega}$,
for some measurable subset $\Omega$ of $\R^n$. The wavelet property is then
closely related to geometric properties of the set $\Omega$, in particular
spectral and tiling properties of $\Omega$. The study of wavelet sets then becomes
an interplay between group theory, geometry, operator theory and analysis,
cf \cite{dai_larson_speegle:97,dai_diao_gu_han:03}.
Our main results are existence theorems for such wavelets for some special
dilation sets $\mathcal {D}$ which are \textit{not necessarily} groups,
see Theorem \ref{th-2.1} and Theorem \ref{th-2.2}.

The article is organized as follows. In Section \ref{s:background}
we recall some well known facts on wavelet sets and their spectral
and tiling properties. In Section 2 we give a construction
of a special type of subspace wavelet set and we continue
with two examples of such sets in Section 3. In Section 4 we discuss
the general framework that has motivated the construction in
this article. Those are the prehomogeneous vector spaces.

\section{Wavelet sets}\label{s:background}
\noindent
In this section we recall some basic facts about wavelet sets.
In particular we discuss their
relation to  \textit{spectral sets}  and to \textit{tilings}.
A standard reference is the work by Y. Wang \cite{wang:02}. See
also the discussion in \cite{olafsson_speegle:03} and the reference therein.

Set $e_{\lambda}(\xi)=e^{2\pi i\langle{\lambda,\xi}\rangle}$. If $\Omega$ is a
Lebesgue measurable
subset of $\R^n$, then
$$|\Omega |=\int \chi_\Omega (t_1,\ldots ,t_n)\, dt_1\ldots dt_n$$
denotes the
measure of $\Omega$ with respect
to the standard Lebesgue measure on $\R^n$.

\begin{de}
A set $\Omega\subseteq \mathbb{R}^n$ with $0< |\Omega| <\infty$ is
a spectral set if there exists a
set $\mathcal{T}\subseteq\mathbb{R}^n $ such that
$\{e_{\lambda}\mid\lambda\in\mathcal{T}\}$
is an orthogonal basis for $L^{2}(\Omega)$.
The set $\mathcal{T}$ is called a spectrum
of $\Omega$, and $(\Omega,\mathcal{T})$ is said to be a spectral pair.
\end{de}

\begin{de}
A measurable tiling of a measure space $(M,\mu)$ is a countable collection of
subsets $\{\Omega_j\}$ of $M$, such that
\begin{displaymath}
\mu(\Omega_i\cap\Omega_j)=0\, ,
\end{displaymath}
for $i\not= j$,
and
$$
\mu(M\diagdown\bigcup_{j}\Omega_j)=0\, .$$
\end {de}
\begin{de} Let $M\subseteq \R^n$ be a measurable set with $| M |>0$.
Let $\mathcal {D}\subseteq \GL (n,\mathbb{R})$ and $\mathcal{T}\subseteq\mathbb{R}^n$.
\begin{enumerate}
\item[1)]
We call
$\mathcal{D}$ a  multiplicative tiling set of $M$ if there exists a
measurable set
$\Omega\subseteq \R^n$,
$|\Omega |>0$, such that $\{d\Omega\mid d\in\mathcal{D}\}$ is a
measurable tiling of $M$.
The set $\Omega $ is called a  multiplicative $\cD$-tile
for $M$
or simply a \textit{multiplicative tile};
\item[2)]
We call  $\mathcal{T}$ an additive tiling set of $\mathbb{R}^n$ if there
exists a measurable set $\Omega\subseteq\mathbb{R}^n$, $|\Omega | >0$, such that
$\{\Omega+t\mid t\in\mathcal{T}\}$ is a measurable tiling of $\mathbb{R}^n$.
The set $\Omega $ is called an \textit{additive $\mathcal{T}$-tile}
or simply an \textit{additive tile};
\item[3)]
A set $\Omega $ is called a $(\mathcal D,\mathcal T)$-\textit{tile}
for $M$ if $\Omega$ is a
$\mathcal{D}$-multiplicative tile for $M$ and a $\mathcal{T}$-additive tiling
set for $\R^n$.
\end{enumerate}
\end{de}
\begin{rem} If $M=\R^n$ then we will simply speak of multiplicative tiling
and $(\cD, \cT)$-tile, without the reference to $M$.
\end{rem}
\begin{rem}\label{re-1} Note, that we do not assume that $\cD M\subseteq M$
or even
\begin{equation}\label{eq-1}
|(M\setminus \cD M)\cup ((\cD M)\setminus M)|=0\, .
\end{equation}
Neither do we assume that there exists a
zero set $Z$, such that $\Omega\subseteq M\cup Z$. That will
always be the case if $\id \in\cD$. It should be noted, that we can always,
without loss of generality, assume that $\id \in \cD$ and
$\Omega\subseteq M$. To see that, fix a $d\in\cD$ and
replace $\cD$ by $\cD d^{-1}$ and $\Omega$ by
$d\Omega\subseteq M$ (up to a set of measure zero).
Note, that if $\cD$ is a group, then $\Omega\subseteq M$,
and (\ref{eq-1}) follows from the tiling property (1).
\end{rem}
The spectral property of a set is closely related to the tiling property,
in particular if the spectrum is a lattice. This was first
noticed by  Fuglede in  \cite{fuglede:74}. For
a non-empty subset $\cT\subset \R^n$, set
$$\cT^\ast :=\{ t\in\mathbb{R}^n \mid \langle{t,s}\rangle
\in \mathbb{Z}, \hbox{ for all } {s}
\in\mathcal{T}\}\, . $$
If $\cT$ is a lattice, then so is $\cT^\ast$. In that case $\cT^\ast$
is called \textit{the dual lattice} of $\cT$.

\begin{thm}[Fuglede \cite{fuglede:74}]
Assume that $\mathcal{T}$ is a lattice. Then $\Omega $ is a spectral set with spectrum $\mathcal{T}$ if
and only if $\{\Omega+t\mid t\in\mathcal{T}^{\ast}\}$ is a measurable tiling of $\mathbb{R}^n$.
\end{thm}
This result and several examples led Fuglede to conjecture, cf. \cite{fuglede:74}:
\begin{conj}[The Spectral-Tile Conjecture]
A measurable set $\Omega$, with positive and finite measure, is a spectral set if and only if it is
an additive tile.
\end{conj}

Several people worked on this conjecture and derived important results and validated the
conjecture for some special cases, see \cite{IKT99,JP98,JP99,LW97,wang:02}
and the references therein.
However, in 2003, T. Tao \cite{tao:04} showed that the conjecture is false in dimension
$5$ and higher if the lattice hypothesis is
dropped. The other direction was disproved by
M. N. Kolountzakis and M. Matolcsi
in \cite{KM04,KM04a}.
But even now, after  the Spectral-Tiling conjecture has
been proven to fail in higher dimensions, it is still interesting and
important to understand better the connection between spectral properties
and tiling in particular, because of the connection to
wavelet sets.

\begin{thm}[Tao \cite{tao:04}] Let $n\ge5$ be an integer. Then there exists a
compact set $\Omega\subset \mathbb{R}^n$
of positive measure such that $L^2(\Omega)$ admits an orthonormal basis of exponentials
$\{e_{\lambda}\mid\lambda\in\mathcal{T}\}$ for some $\mathcal{T}\subset\mathbb{R}^n$, but
such that $\Omega $  can not tile $\R^n$ by
translation.  In particular,
Fuglede's conjecture is false in $\mathbb{R}^n$ for $n\ge5$.
\end{thm}

We now discuss briefly some results by Dai, Larson, and Speegle \cite{dai_larson_speegle:97}.
For that assume that $\mathcal{G}$  and $\mathcal{T}$
are countable groups acting on a measure space $(M,\mu)$ by
measurable automorphisms. We think of $\mathcal{G}$ as
``dilation'' group and $\mathcal{T}$ as
translations. We write correspondingly the action
of $\mathcal{G}$ by $m\mapsto gm$ and the action  of
$\mathcal{T}$ by $m \mapsto t+m$.
Two sets $E$ and $F$ are said to be $\mathcal{G}$-\textit{dilation congruent},
$E\sim_{\mathcal{G}}F$, if there
exist measurable partitions $\{E_i\}$ and  $\{F_i\}$ of $E$ and $F$, respectively, such that
$F_i=g_iE_i$ for some  $g_i\in\mathcal{G}$. Similarly,
two sets $E$ and $F$ are said to be $\mathcal{T}$-\textit{translation congruent},
$E\sim_{\mathcal{T}}F$, if there
exist measurable partitions $\{E_i\}$ and
$\{F_i\}$ of $E$ and $F$, respectively, such that
$F_i=t_i+E_i$ for some  $t_i\in\mathcal{T}$.

In 1996, the above three authors showed in \cite{dai_larson_speegle:97},
that wavelet sets exist for  groups of dilations. They first introduce the notion of
\textit{abstract dilation-translation pair}.

\begin{de}[Dai-Larson-Speegle \cite{dai_larson_speegle:98}]
Let $M$ be a metric space and
$\cD$ and $\cT$ discrete groups of automorphisms of $M$.
A pair $(\mathcal{D},\mathcal{T})$ is called an abstract
dilation-translation pair if the following holds:

\begin{enumerate}
\item[1)]For each bounded set $E$ and each open set
$F$ there exist $d\in \mathcal{D}$ and $t\in \mathcal{T}$
such that $ t + E\subseteq d(F)$;
\item[2)] There is a fixed point $\theta $ for $\mathcal{D}$
such that for any neighborhood $N$
of $\theta $ and for any bounded set $E$, there is an element
$d\in \mathcal{D}$ such that $d(E)\subseteq N $.
\end{enumerate}
\end{de}
In \cite{dai_larson_speegle:97} the following is proved:
\begin{thm}[Dai-Larson-Speegle \cite{dai_larson_speegle:97}]\label{th-dls1}
Let $M$ be a metric space and
$(\mathcal{D},\mathcal{T})$  an abstract
dilation-translation pair with $\theta$ a $\mathcal{D}$ fixed point
as above. If $E$ and $F$ are bounded measurable sets in $M$ such that
$E$ contains a neighborhood of $\theta$
and $F$ has nonempty interior and is bounded away
from $\theta$, then there exists a measurable set
$W\subseteq M$, $W\subseteq\cup_{d\in\mathcal{D}}d(F)$
which is $\mathcal{D}$-dilation congruent to $F$ and
$\mathcal{T}$-translation congruent to $E$.
\end{thm}
We now apply this abstract result to $\R^n$.
If $d\in \GL (n,\R)$, $\gamma\in \R^n$, and $\psi :\R^n\to \C$
is measurable, set
\begin{equation}\label{eq-psidg}
\psi_{d,\gamma}(x)=|\det d|^{1/2}\psi(dx+\gamma)\, .
\end{equation}
Note that the
Fourier transform of $\psi_{d,\gamma}$ is given by
\begin{equation}\label{eq-ft}
\widehat{\psi_{d,\gamma}}(\lambda )=e^{2\pi i <\gamma ,d^{-T}\lambda>}
\widehat{\psi} (d^{-T}\lambda )\, .
\end{equation}

\begin{de} Let $M\subseteq \R^n$ be measurable, $|M|>0$, and
$\cD\subset \GL (n,\R)$. Let $\cT\subset \R^n$ be
discrete. Then a measurable set $\Omega\subseteq M$
is called an $M$-subspace $(\cD,\cT)$-wavelet set if the set
of functions
$\{\psi_{d,\gamma}\}_{(d,\gamma)\in\cD\times \cT}$ is an orthogonal
basis for $L^2_M(\R^n)$, where $\psi =\cF^{-1}\chi_\Omega$.
\end{de}

\begin{rem}\label{re-2} Note again that we do not assume that $\Omega\subset M$, nor
that $\cD M =M$ up to set of measure zero, but this follows if
$\id \in \cD$. As in Remark \ref{re-1} one can always assume this
by replacing $\cD$ by $d^{-1}\cD$ and $\Omega$ by $d^{T}\Omega$ for
a fixed $d\in \cD$.
\end{rem}

We get from Theorem \ref{th-dls1}:

\begin{thm}[Dai-Larson-Speegle \cite{dai_larson_speegle:97}]\label{th-DLS}
Let $a$ be an expansive matrix, and let $M\subseteq \mathbb{R}^n$ be a measurable
set of positive
measure such that $a^{T}M=M$. Let $\mathcal{D}=\{a^k\mid k\in\mathbb{Z}\}$ and
let $\mathcal{T}$ be a
full rank lattice. Then there exists an $M$-subspace $(\cD,\cT)$-wavelet set.
\end{thm}

There are several generalizations of this Theorem. We
refer to \cite{olafsson_speegle:03}
for discussion and references. We will only mention two
important results here:

\begin{thm}[Dai-Diao-Gu-Han \cite{dai_diao_gu_han:03}]\label{thm8}
Let $M$ be a measurable subset of $\mathbb{R}^n$, with positive measure satisfying
$a^TM=M$, for
some expansive matrix $a$, and let
$\mathcal{T}$ be a full rank lattice. Then there exists a set $\Omega
\subseteq M$ such that
$ \{\Omega + t\mid t\in \mathcal{T}\}$ is a measurable tiling of $\mathbb{R}^{n}$
and $\{ (a^T)^{k}\Omega \mid k\in \mathbb{Z}\}$  is a measurable tiling of $M$. In
particular $\Omega$ is an $M$-subspace $(\cD,\cT)$-wavelet set.
\end{thm}

\begin{thm}[Wang \cite{wang:02}]\label{th-wang}
Let $\mathcal {D}\subseteq GL(n,\mathbb{R})$ and $\mathcal{T}\subseteq\mathbb{R}^n $. Let
$\Omega\subseteq\mathbb{R}^n$ be measurable, with positive and finite measure.
 If $\Omega$ is a measurable $\mathcal{D}^T$-tile and $(\Omega,\mathcal{T})$ is a spectral pair,
 then $\Omega $ is a $(\mathcal D,\mathcal T)$-wavelet set.
Conversely, if $\Omega $ is a $(\mathcal D,\mathcal T)$-wavelet set and $0\in\mathcal T$, then
$\Omega$ is a measurable $\mathcal{D}^T$-tile and $(\Omega,\mathcal{T})$ is a spectral pair.
\end{thm}

Let us sketch some of the ideas of the proof to underline the
connection between spectral properties, tiling, and wavelet sets.
Let $\psi=\mathcal{F}^{-1}\chi_{\Omega}$.
As the Fourier transform is an unitary isomorphism,
it follows, that the
set $\{\psi_{d,t}\mid d\in\mathcal{D},t\in\mathcal{T}\}$ is an orthogonal basis for
$L^{2}_M(\mathbb{R}^n)$ if and only if
the set $\{\widehat{\psi_{d,t}}\mid d\in\mathcal{D},t\in\mathcal{T}\}$ is
an orthogonal basis for $L^2(M )$. Here, as before, we have set
$$\psi_{d,t}(x)=|\det d|^{1/2}\psi(dx+t)\, .$$
A simple calculation shows that
\begin{displaymath}
\widehat{\psi_{d,t}}(\lambda )
={|\det d|}^{-\frac{1}{2}}e^{2\pi i\langle{d^{-1}t,\lambda }\rangle}\chi_{d^{T}\Omega}(\lambda)
={|\det d|}^{-\frac{1}{2}}e^{2\pi i\langle{t,d^{-T}\lambda }\rangle}\chi_{d^{T}\Omega}(\lambda)\, .
\end{displaymath}
The fact, that $d^T\Omega$ is a measurable tiling of $M$ implies
that
$$L^2(M )\simeq \bigoplus_{d\in \cD}L^2(d^T\Omega )\, .$$
The orthogonal projection onto $L^2(d^T\Omega)$ is given by
$f\mapsto f\chi_{d^T\Omega}$ and
$$f=\sum_{d\in \cD} f\chi_{d^T\Omega }\, .$$
The spectral property implies, that
$\{e_t\}_{t\in\cT}$ is an orthogonal basis for $L^{2}(\Omega)$.
As the linear map
$f\mapsto |\det d |^{-1/2}f(d^{-T}\cdot )$ is a unitary
isomorphism $L^2(\Omega )\simeq L^2(d^T \Omega)$ it follows, that
the set of functions $\{|\det d |^{-1/2} e^{2\pi i\langle{t,d^{-T} \cdot }\rangle}
=
|\det d|^{-1/2}e_{d^{-1} t}
\mid t\in\mathcal{T}\}$ is an orthogonal
basis for $L^{2}(d^T\Omega )$. Putting those two things together, implies
that $ \{|\det d |^{-1/2} e_{d^{-1}t}
\mid d\in\cD,\,  t\in\mathcal{T}\}$ is an orthogonal basis for
$L^2(M)$. Hence $\{\psi_{d,t}\mid (d,t)\in \mathcal{D}\times \mathcal{T}\}$
is an orthogonal basis for $L^2_M(\R^n)$.

We refer to \cite{olafsson_speegle:03}
for further discussion and generalizations of Theorem \ref{th-wang}.

\section{Existence of subspace wavelet sets}\label{s-2}
\noindent
In this section we discuss how to construct subspace wavelet sets using a
kind of ``induction'' process, i.e., using well known facts discussed in
the previous section on smaller dilation sets acting on a smaller
frequency set and then extending those to our bigger dilation set and
frequency set.
We start with two simple, but important, observations.
For $\cA,\cB\subset \GL (n,\R)$ we say that the product
$\cA\cB=\{ab\mid a\in \cA,\, b\in\cB\}$ is \textit{direct} if
$a_1b_1=a_2b_2$, $a_1,a_2\in\cA$, $b_1,b_2\in\cB$, implies that
$a_1=a_2$ and $b_1=b_2$.

We state the following simple Lemma, but note, that we will be
using the proof more than the actual statement. We remark, that
the statements in the next few lemmas and theorems hold for
more general settings, i.e., one could replace $\R^n$ by a measure space
and $\GL (n,\R)$ by a group of automorphisms of $M$.

\begin{lem}\label{th-2.1}
Let $M\subseteq \R^n$ be measurable.
Let $\mathcal{A},\mathcal{B}\subset \GL(n,\mathbb{R})$ be two non-empty sets, such that
the product $\cA\cB$ is direct
Let $\mathcal{D} = \mathcal{A}\mathcal{B}=\{ab\mid a\in
\mathcal{A}\, ,\,\, b\in \mathcal{B}\}$.
Then there exists a $\mathcal{D}$-tile $\Omega$ for $M$ if and only if there
exists a measurable set $N\subseteq \R^n$,
such that $\mathcal{A} N$ is a measurable tiling of
$M$, and a
$\mathcal{B}$-tile $\Omega$ for $N$.
\end{lem}
\begin{proof}
Assume that $\Omega$ is a $\cD$-tile for $M$. Set
$$N:= \mathcal{B}\Omega=\bigcup_{b\in\cB}b\Omega\, .$$
Assume, that there are $b_1,b_2\in\cB$ such
that $|b_1\Omega \cap b_2\Omega|>0$.
Then $|(ab_1\Omega)\,  \cap \, (ab_2\Omega)|>0$ for all $a\in\cA$, which
contradicts our assumption, that $\cD \Omega$ is a measurable
tiling of $M$. Hence $\cB\Omega$ is a measurable tiling of $N$.
We have, up to a set of measure zero:
$$\cA N=\bigcup_{a\in \cA}\, aN =\bigcup_{a\in\cA ,\, b\in\cB}\, ab\Omega
=M\, .$$
Assume, that there are $a_1,a_2\in\cA$ such that $|a_1 N\cap a_2 N|>0$.
Then we can find $b_1,b_2\in \cB$ such that
$|a_1b_1\Omega\cap a_2b_2\Omega|>0$. As the product $\cA\cB$ is direct,
and $\cD \Omega$ is a measurable tiling of $M$, it follows that
$a_1b_1=a_2b_2$ and hence $a_1=a_2$.
Thus $\cA N$ is a measurable tiling of $M$.

For the other direction, assume that $\cA N$ is
a measurable tiling of $M$ and that $\cB \Omega$ is a measurable
tiling of $N$.
Then (up to sets of measure zero)
$$\bigcup_{d\in \cD}\, d\Omega=\bigcup_{a\in\cA}\bigcup_{b\in \cB}ab\Omega
=\bigcup_{a\in\cA}a \bigcup_{b\in \cB}b\Omega
=\bigcup_{a\in\cA}aN = M\, .$$
Assume that $|d_1\Omega\cap d_2\Omega|>0$. Then there
are unique $a_1,a_2\in\cA$, and $b_1,b_2\in \cB$ such that
$d_1=a_1b_1$ and $d_2=a_2b_2$. Hence
$|a_1N\cap a_2N|>0$ which implies  that $a_1=a_2$, as $\cA N$ is a
measurable tiling for $M$.
But then $|b_1\Omega \cap b_2\Omega|>0$, which implies that $b_1=b_2$. Hence
$d_1=d_2$. This shows, that $\cD \Omega$ is a measurable tiling
of $M$.
\end{proof}

\begin{rem} We would like to remark at this point, that we do
not assume that $\Omega\subseteq M$, nor that $N\subseteq M$.
But this will in fact be the case in most applications because
$\cD$ will contain the identity matrix. Recall also
from Remark \ref{re-1} and Remark \ref{re-2} that
we can always assume that $\id\in\cD$ and $\Omega\subseteq M$
up to set of measure zero. The same remarks hold for
the following Theorems.
\end{rem}

\begin{thm}[Construction of wavelet sets by steps, I]\label{th-2.2}
Let $\mathcal{M},\mathcal{N}\subset \GL(n,\R)$ be two non-empty
subsets such that the product $\cM\cN$ is direct.
Let $\mathcal{L} = \mathcal{M}\mathcal{N}$.
Assume that $M\subseteq \mathbb{R}^n$ with $|M|>0$, is measurable.
Let $\mathcal{T}\subset \mathbb{R}^n$ be discrete.
Then there exists a $(\mathcal{L} ,\mathcal{T})$-wavelet
set $\Omega\subset M$ for $M$ if and only if there
exists a $\mathcal{N}^T$-tiling set $N\subset M$ and a
$(\mathcal{M},\mathcal{T})$-wavelet set $\Omega_1$ for $N$.
\end{thm}
\begin{proof}
Set $\cA =\cN^T$ and $\cB=\cM^T$.
Then the conditions in Lemma \ref{th-2.1} are satisfied.

Assume that $\Omega\subset M$ is a $(\mathcal{L}, \mathcal{T})$-wavelet
set for $M$. As above we set
$$N:=\cB \Omega:=\bigcup_{b\in\mathcal{M}}\, b^T \Omega\, .$$
Then, as above,  we see that
$\cA N$ is a measurable tiling of $M$.
As  $\Omega$ is a spectral set,
it follows from Theorem \ref{th-wang} that $\Omega$
is a $(\cM,\cT)$-wavelet set for $N$.

Assume now that $N$ is a $\mathcal{N}^T$-tiling for $M$,
and that $\Omega_1$ is a $(\mathcal{M},\mathcal{T})$-wavelet set for $N$. Then,
in particular $\Omega_1$ is a $\cB$-tile for $N$. As
$\cA N$ is a measurable tiling of $M$, it follows that $\cL^T\Omega_1$ is
a measurable tiling of $M$. As $\Omega_1$ is a spectral set
it follows from Theorem \ref{th-wang} that
$\Omega_1$ is a $(\mathcal{L} ,\mathcal{T})$-wavelet set for $M$.
\end{proof}
Recall that if $\cD\subseteq \GL(n,\R)$, and $\cG\subset \GL (n,\R)$ is
a group that acts on $\cD$ from the right, then there exists
a subset $\cD_1\subseteq \cD$, such that $\cD=\cD_1\cG$ and the
product is direct. Note, that we do not assume that $\cG\subset \cD$.

\begin{thm}[Construction of wavelet sets by steps, II] Let $\cD\subset
\GL (n,\R)$ and $M\subseteq \R^n$ measurable
with $|M|>0$. Let $\cT\subset \R^n$ be discrete.
Assume that $\mathcal{G}\subset \GL(n,\R)$ is a group that
acts on $\cD$ from the right. Let $\cD_1\subseteq \cD$ be
such that $\cD=\cD_1\cG$ as a direct product.
Then there exists a $(\cD,\cT)$-wavelet set
$\Omega$ for $M$ if only only if there exists
a $\cG^T$-tiling set $N$ for $M$ and a $(\cD_1,\cT)$-wavelet set $\Omega_1$ for
$N$.
\end{thm}
\begin{proof} This follows directly from Theorem \ref{th-2.2} with
$\cM=\cD_1$ and $\cN=\cG$.
\end{proof}

The question is then, how to obtain a wavelet set for
the starting subset $N$. The following gives one
way to do that.

\begin{lem}[Existence of subspace wavelet sets]\label{th-exwave}
Let $M \subseteq \mathbb{R}^{n}$ be a measurable
set, $|M|>0$. Let $a\in \GL(n,\R)$ be an expansive
matrix and $\emptyset\not=\mathcal{D}\subset \GL (n,\mathbb{R}) $.
Assume that
$\cD^T$ is  a multiplicative tiling of $M$, $a\cD=\cD$ and
$a^T M=M$. If $\mathcal{T}$ is a lattice,
then there exists a
measurable set $\Omega \subseteq M$ such that $ \Omega + \cT$ is a
measurable tiling
of $\mathbb{R}^{n}$
and $\cD^T \Omega$  is a measurable tiling of
$M$. In particular, $\Omega$ is a $M$-subspace $(\cD,\cT)$-wavelet set.
\end{lem}

\begin{proof} Let $b=a^T$ and $\cB=\{b^k\mid k\in\Z\}$. Then
$\cB$ is an abelian group that acts on $\cD^T$
from the right. Hence,
there exists a set $\cA\subset \cD$ such that
$\cD^T=\cA\cB$ and the product is direct. Thus,
the conditions in the previous Theorems are satisfied.

Let $E\subset M$ be such that $\cD^T E$ is a measurable
tiling of $M$. Set $N:=\cB E\subseteq M$. Then
$N$ is $\cB$ invariant. By Theorem \ref{th-DLS} there
exists a $(\cB,\cT)$-wavelet set $\Omega$ for $N$.
Set
$$N:=\bigcup_{k\in \Z}\, b^k\Omega\, .$$
Then, as before, we see that $\Omega$ is a $(\cD,\cT)$-wavelet.
\end{proof}

\section{Two simple examples}
\noindent
Before we discuss more general situations in
the next section, let
us explain our results by two simple
examples.
\begin{Example}\label{Ex-3.1}
For $\theta\in \R$, let $R_\theta=\begin{pmatrix} \cos \theta & -\sin \theta
\cr \sin \theta &\cos \theta\end{pmatrix}$
denote the rotation in $\R^2$ by the angle $\theta$.
Let $a>1$. For any integer $m\ge 2$ let
 \begin{displaymath}
\mathcal{D}_{a,m}:=\{a^{n}R^{k}_{2\pi /m}|n\in\mathbb{Z},\ k=0,\dots,m-1\}
\, .
\end{displaymath}
Note, that $R^k_{2\pi /m}=R_{2\pi k/m}$ and that $\cD_{a,m}$
is a group.
Let
$\mathcal{T}=\mathbb{Z}^2$
and
$$\R^2_{2\pi/m}=\{r (\cos \psi,\sin \psi)^T\mid 0\le \psi \le {2\pi/m}\, ,\,\,
r>0\}\, .$$
Then $\R^2_{2\pi/m}$ is a tiling set for the finite group
$\{R^k_{2\pi /m}\mid k=0,1,\ldots ,m-1\}$.
As $a\cdot\id$ is expansive and $a\R^2_{2\pi /m}=\R^2_{2\pi /m}$,
it follows from Theorem \ref{th-DLS}
that there exists a $(\cA:=\{a^k
\cdot \id \mid k\in\Z\},\cT)$-wavelet
set $\Omega$ for $\R^2_{2\pi/m}$ and hence a $(\cD_{a,m},\cT)$-wavelet
set for $\R^2$, see also  Lemma \ref{th-exwave}.
We show here how to construct such a wavelet set. Note that
we only have to construct a $(\cA,\cT)$-wavelet set for
$\R^2_{2\pi/m}$.
For that,
let
$$E =[0,1]\times [0,\tan({2\pi/m})]\,\,\, \hbox{if} \,\, m\neq2,4$$
$$E =[0,1]^2\,\,\, \hbox{if} \,\, m=4$$
$$E =[-1,1]\times[0,1]\,\,\, \hbox{if} \,\, m=2$$
$$F =\{(x,y)\in \R^2_{2\pi/m}| 1<x<a\}.$$

It is clear that $(\mathcal{D}_{a,m},\mathcal{T})$ is an abstract
dilation-translation pair with a fixed point $0$.
By  Theorem 1 in \cite{dai_larson_speegle:97}, it
follows that there exists a measurable set
$W$ such that $W$ and $E$ are $\mathcal{T}$-\textit{translation congruent},
and $W$ and $F$ are $\mathcal{D}$-\textit{dilation congruent}.
 On the other hand, $F$ is a
$\mathcal{D}_{a,m}$-multiplicative tile and $\{E,\mathcal{T} \}$ is a spectral pair.
It follows that $W$ is a $\mathcal{D}_{a,m}$-multiplicative tile
and $\{W,\mathcal{T}\}$ is a spectral set.
Thus, by Theorem \ref{th-wang}  $W$ is a $(\mathcal{D}_{a,m}, \mathcal{T})$
wavelet set.

\begin{figure}[!ht]
\begin{center}
\epsfig{scale=0.7,angle=0,file=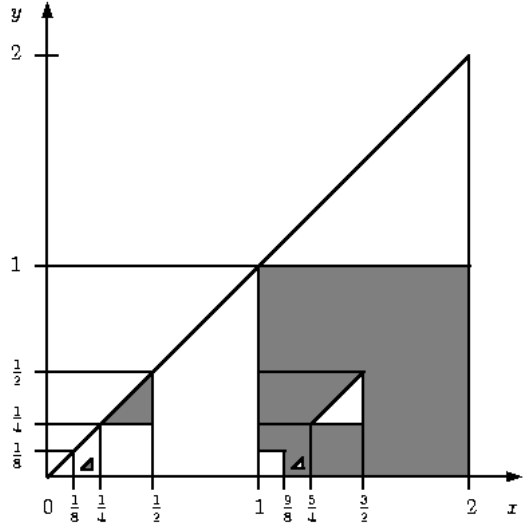}
\caption{A $(\mathcal{D}_{2,4},\mathcal{T})$ wavelet set.}
\label{fig1}
\end{center}
\end{figure}
.
\end{Example}

\begin{Example}
Let
$A=(\begin{smallmatrix}
2&0\\0&3
\end{smallmatrix})$, $\mathcal{D}=\{R^{k}_{{\pi}/{2}}A^n\mid n\in\mathbb{Z},
k=0,1,2,3\}$ and $\mathcal{T}=\mathbb{Z}^2$. Note that $\mathcal{D}$
is not a group anylonger. Let $E=[0,1]^2$ and $F=([0,2]\times[0,3])\setminus[0,1]^2$. Then $E$ is a $\mathbb{Z}^2$-additive tile
and that $F$ is a $\mathcal{D}$-multiplicative tile. Following the same procedure,
we get a set $W$ such that $W$ and $E$ are $\mathbb{Z}^2$-\textit{translation congruent}, $W$ and $F$ are $\mathcal{D}$-\textit{dilation congruent} and so, it follows that $W$ is a $\mathcal{D}$-multiplicative tile
and $\{W,\mathbb{Z}^2\}$ is a spectral set.
Thus, $W$ is a $(\mathcal{D}, \mathbb{Z}^2)$
wavelet set.
The wavelet set $W$ has the form
\begin{displaymath}
W=\bigcup_{i=1}^{2}\bigcup_{j=1}^{\infty}W_{i,j}\, .
\end{displaymath}
The description of the set $W_{i,j}$ is as follows
\begin{displaymath}
\begin{array}{rcl}
W_{1,1} & = & (E\setminus A^{-1}E)+(1,0)\\[1ex]
W_{2,1} & = & A^{-2}[(0,1)\times(1,3)]\\[1ex]
\end{array}
\end{displaymath}
For $j\geq 2$, we have the following formulas
\begin{displaymath}
W_{1,j}=[(A^{-j+1}E\setminus A^{-j}E)\setminus W_{2,j-1}]+(1,0)
\end{displaymath}
\begin{displaymath}
W_{2,j}=A^{-j-1}[W_{2,j-1}+(0,1)].
\end{displaymath}

\begin{figure}[!ht]
\begin{center}
\epsfig{scale=0.7,angle=0,file=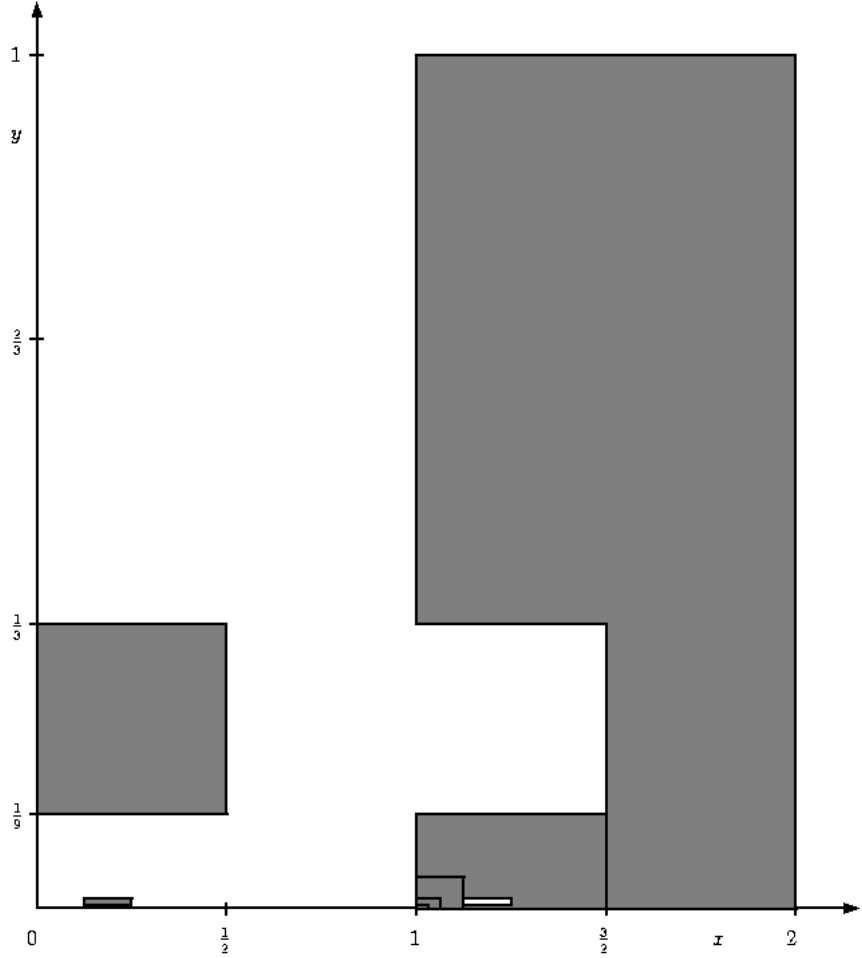}
\caption{The set $W$.}
\label{fig3}
\end{center}
\end{figure}
\end{Example}
\section{More general examples}
\noindent
The results in Section \ref{s-2} give us tools to
construct wavelet sets or subspace wavelet sets.
In this section we discuss some general settings for
our results. Those examples were discussed in
\cite{FO02,O02,olafsson_speegle:03} and the
results therein were one of the main motivations
for starting this work. We refer to the above
articles for more information and details.  The
starting point is a subgroup $H\subseteq \GL (n,\R)$
such that $H^T$ has finitely many open orbits $U_1,\ldots , U_k$
in $\R^n$ such that
$$|\R^n\setminus (U_1\cup \ldots \cup U_k)|=0\, .$$
In this case, the pair $(H,\R^n)$ is  called a
\textit{Prehomogeneous vector
space}.
Furthermore, we assume that
$$H^u=\{h\in H\mid h^T(u)=u\}$$
is compact for all $u\in U_j$, $j=1,\ldots ,k$.
The reason for those assumptions is, that this
implies that $H$ is admissible in the sense of \cite{LWWW2002,WW01}, and
we can define the continuous wavelet transform. In fact, there exists a
function $\psi \in L^2(\R^n)$ such that
$$\int_H |\psi (h^T x)|^2\, dh=1$$
for almost all $x\in \R^n$. Here $dh$ stands for a left-invariant
Haar measure on $H$. Note, that we need the transpose of $h$ here,
because of the action in the frequency domain, cf (\ref{eq-ft}).

Let $G:=\R^n\times_s H$ be the
semi-direct product of
the groups $\R^n$ and $H$. This is a natural generalization
of the $(ax+b)$-group in one dimension. Define
$W_\psi : L^2(\R^n)\to L^2(\R^n\times_s H,\frac{dydh}{|\det h|})$ by
\begin{equation}\label{eq-inversion}
W_\psi (f)(y,h)=|\det h|^{-1/2}\int_{\R^n} f(x)\overline{\psi (h^{-1}(x-y))}\, dx\, .
\end{equation}
That $W_\psi: L^2(\R^n)\to L^2(G)$ is an unitary isomorphism onto
its image in $L^2(G)$ was shown in \cite{FO02}. It is the
\textit{continuous wavelet transform}. The
inverse of $W_\psi$ is given by the weak integral:
\begin{equation}\label{eq-inv}
f=\int_{G}W_\psi (f)(y,h) |\det h|^{-1/2}\psi (h^{-1}(\cdot -y)\, \frac{dydh}{|\det h|}\, .
\end{equation}
We refer to \cite{FO02} and Section 2 of \cite{olafsson_speegle:03} for
detailed discussion and proofs.
The natural question now arises to find a discrete version of
the inversion formula (\ref{eq-inv}). Note, that we can replace
$L^2(\R^n)$ by the subspaces $L^2_{U_j}(\R^n)$ and define a
continuous wavelet transform on each of those spaces. This is where the
connection to the first part of this article comes into play.
\begin{Example} Let $n=1$ and $H=\R^+$ acting on $\R$ by
multiplication. Then $G$ is just the $(ax+b)$-group. We have
two open $H$-orbits $\R^+$ and $\R^-:= -\R^+$. In each case the
stabilizer of a point is trivial, i.e. just the one point
subgroup $\{1\}$ and hence compact. In this case the
space $L^2_{\R^+}(\R )$ is the classical Hardy space of all
functions in $L^2(\R )$ that can be realized as $L^2$-boundary values of
holomorphic functions on the upper half-plane $\R + i\R^+$. Similarly,
the space $L^2_{\R^-}(\R )$ is the space of all $L^2$-functions on the
line, which can be written as $L^2$-boundary values of holomorphic
functions on the lower half-plane $\R+i\R^-$.
If we replace $\R^+$ by $\R^*=\R\setminus \{0\}$, then there is
only one open orbit $\R^*$.
\end{Example}
\begin{Example} Let now $H=\R^+\SO (n)$. Then we have only one
open orbit $\R^n\setminus \{0\}$. The stabilizer of a point
in $\R^n\setminus \{0\}$ is naturally isomorphic to $\SO (n-1)$ and
hence compact. A discrete version of
the wavelet transform in this case would be a generalization of
Example \ref{Ex-3.1}.
\end{Example}

\begin{Example} Let $\Omega\subset \R^n$ be an open convex cone.
Define
$$G(\Omega ):= \{h\in \GL (n,\R)\mid h(\Omega )=\Omega\}$$
and
$$\Omega^*:=\{w\in \R^n \mid \forall v\in \overline{\Omega}\setminus \{0\}\, :\,
\langle{w,v}\rangle >0\}\, .$$
The cone $\Omega$ is said to be homogeneous if $G(\Omega )$ acts
transitively on $\Omega$ and self-dual if $\Omega^*=\Omega$.
The cone $\Omega$ is symmetric if it is homogeneous and self-dual.
Then there exists a solvable Lie group $H=AN\subset \GL (n,\R )$,
such that $H$ acts freely and transitively on $\Omega$.
Here $N$ is simply connected and nilpotent, the group $A$ is
simply connected and abelian, and $A$ normalizes $N$. Note
that $\R^+\id\subset A$ is a one dimensional subgroup containing
an expansive matrix. Furthermore $\R^+\id$ commutes with all elements
in $H$. This example is important in the study
of Hardy spaces on $\R^n$ as well as Besov-spaces
associated to symmetric cones, c.f.,
\cite{BBGR}
\end{Example}

As in Section 3 of \cite{olafsson_speegle:03} we consider
the case where the group  $H$ can be written as a
semi direct product $H=ANR=NAR=RAN$ such that
\begin{enumerate}
\item[1)] $R$ is
compact, commutes with $A$, and
normalizes $N$;
\item[2)] $A$ is simply connected abelian and normalizing
$N$;
\item[3)] the
map
\begin{equation*}
N\times A\times R\ni (n,a,r)\mapsto nar\in H
\end{equation*}
is a diffeomorphism;
\item[4)] There exists an expansive matrix $a\in A$ which is central in
the group $H$;
\item[5)] there exists a co-compact discrete subgroup $\Gamma _{N}\subset
N$.
\end{enumerate}
In examples one of the groups can be trivial, and one could also
use the group $AR$ to fit the examples directly to
the previous sections. Those assumptions includes  all the examples above.

Assume now that $M$ is one of the open $H^T$-orbits. For simplicity
we assume that $H^T$ acts freely on $M$, otherwise we will have
to be a little more careful about the choice of $\Gamma_R$.
Let $\mathbb{F}_N$ be a compact (or precompact, measurable) subset of $N$, such
that $N=\mathbb{F}_N\Gamma_N$ is a measurable tiling of $N$.
Let $\Gamma _{A}\subset A$ be a
co-compact subgroup in $A$. Choose compact (or precompact measurable)
subsets $\mathbb{F}_{A}\subset A$ such that $A=\mathbb{F}_A\Gamma _{A}$
is a measurable tiling.
Finally let $\Gamma_R\subset R$ be a finite subgroup and
$\mathbb{F}_R$ a compact (or precompact measurable) subset
such that $\mathbb{F}_R\Gamma_R$ is a measurable tiling
of $R$. Set $\mathbb{F}_H:=\F_R\F_A\F_N$ and
$\cD:=\Gamma_N\Gamma_A\Gamma_R$.
Let $\F:=\F_H^Tm_0$, where $m_0\in M$. Then
$\cD^T\F$ is a measurable tiling of $H$ (c.f.
the proof of Lemma 3.2 in \cite{olafsson_speegle:03}).
Let $a$ be the expansive matrix in (4) above.
Then $a\cD =\cD$ and the assumptions in
Lemma \ref{th-exwave} are satisfied. It follows, that
if $\cT$ is any full rank lattice in $\R^n$, then there
exists a $(\cD,\cT)$-wavelet set for
$M$. Using the corresponding wavelet function $\psi=|\Omega |^{-1/2} \cF^{-1}\chi_{\Omega}$
allows us to write a discrete version of the inversion
formula (\ref{eq-inv}):
\begin{eqnarray*}
f &=&\sum_{\gamma=(y,h) \in\cT\times \cD}W_\psi f(\gamma^{-1})\psi_{\gamma^{-1}}\\
&=&\sum_{(y,h)\in\cT\times \cD}|\det h|^{1/2}W_\psi f(-h^{-1}y,h^{-1})\psi ( (\cdot ) + y)\, .
\end{eqnarray*}
The factor $|\Omega |^{-1/2}$ is included to make $\psi$ into an orthonormal
wavelet function.

Note, that the assumption that the central element $a$ as above exist
is not satisfied for arbitrary groups $H$, even if
$H$ has finitely many open orbits, whose union is dense.
As an example take the group
$$H=\left\{\begin{pmatrix} a & t\cr 0 & 1/a\end{pmatrix}
\mid a>0,\, t\in\R\right\}$$
acting on $\R^2$. There are two open orbits, $y>0$, and $y<0$.
Here $A$ is the subgroup of diagonal matrices, and $N$ is
the subgroup of upper triangular matrices with $1$s on the
diagonal. In this case, the center of $H$ is trivial. Note, that $H$
is isomorphic to the $(ax+b)$-group, hence this is in fact
an example of the $(ax+b)$-group acting on $\R^2$.


\begin{thebibliography}{10}

\bibitem{A03}
A.~Aldroubi, C.~Cabrelli, and U.~M. Molter.
\newblock Wavelets on irregular grids with arbitrary dilation matrices and
  frame atoms for {$L\sp 2({\mathbb{R}}\sp d)$}.
\newblock {\em Appl. Comput. Harmon. Anal.}, 17(2):119--140, 2004.

\bibitem{ACDL98}
P.~Aniello, G.~Cassinelli, E.~De~Vito, and A.~Levrero.
\newblock Wavelet transforms and discrete frames associated to semidirect
  products.
\newblock {\em J. Math. Phys.}, 39(8):3965--3973, 1998.

\bibitem{BCMO95}
L.~Baggett, A.~Carey, W.~Moran, and P.~Ohring.
\newblock General existence theorems for orthonormal wavelets, an abstract
  approach.
\newblock {\em Publ. Res. Inst. Math. Sci.}, 31(1):95--111, 1995.

\bibitem{BBGR}
D.~B{\'e}koll{\'e}, A.~Bonami, G.~Garrig{\'o}s, and F.~Ricci.
\newblock Littlewood-{P}aley decompositions related to symmetric cones and
  {B}ergman projections in tube domains.
\newblock {\em Proc. London Math. Soc. (3)}, 89(2):317--360, 2004.

\bibitem{BT96}
D.~Bernier and K.~F. Taylor.
\newblock Wavelets from square-integrable representations.
\newblock {\em SIAM J. Math. Anal.}, 27(2):594--608, 1996.

\bibitem{dai_diao_gu_han:03}
X.~Dai, Y.~Diao, Q.~Gu, and D.~Han.
\newblock The existence of subspace wavelet sets.
\newblock {\em J. Comput. Appl. Math.}, 155(1):83--90, 2003.
\newblock Approximation theory, wavelets and numerical analysis (Chattanooga,
  TN, 2001).

\bibitem{dai_larson:98}
X.~Dai and D.~R. Larson.
\newblock Wandering vectors for unitary systems and orthogonal wavelets.
\newblock {\em Mem. Amer. Math. Soc.}, 134(640):viii+68, 1998.

\bibitem{dai_larson_speegle:97}
X.~Dai, D.~R. Larson, and D.~M. Speegle.
\newblock Wavelet sets in {$\R\sp n$}.
\newblock {\em J. Fourier Anal. Appl.}, 3(4):451--456, 1997.

\bibitem{dai_larson_speegle:98}
X.~Dai, D.~R. Larson, and D.~M. Speegle.
\newblock Wavelet sets in {$\R\sp n$}. {II}.
\newblock In {\em Wavelets, multiwavelets, and their applications (San Diego,
  CA, 1997)}, volume 216 of {\em Contemp. Math.}, pages 15--40. Amer. Math.
  Soc., Providence, RI, 1998.

\bibitem{FO02}
R.~Fabec and G.~{\'O}lafsson.
\newblock The continuous wavelet transform and symmetric spaces.
\newblock {\em Acta Appl. Math.}, 77(1):41--69, 2003.

\bibitem{fuglede:74}
B.~Fuglede.
\newblock Commuting self-adjoint partial differential operators and a group
  theoretic problem.
\newblock {\em J. Functional Analysis}, 16:101--121, 1974.

\bibitem{F96}
H.~F{\"u}hr.
\newblock Wavelet frames and admissibility in higher dimensions.
\newblock {\em J. Math. Phys.}, 37(12):6353--6366, 1996.

\bibitem{F98}
H.~F{\"u}hr.
\newblock Continuous wavelet transforms with abelian dilation groups.
\newblock {\em J. Math. Phys.}, 39(8):3974--3986, 1998.

\bibitem{HW89}
C.~E. Heil and D.~F. Walnut.
\newblock Continuous and discrete wavelet transforms.
\newblock {\em SIAM Rev.}, 31(4):628--666, 1989.

\bibitem{IKT99}
A.~Iosevich, N.~H. Katz, and T.~Tao.
\newblock Convex bodies with a point of curvature do not have {F}ourier bases.
\newblock {\em Amer. J. Math.}, 123(1):115--120, 2001.

\bibitem{JP98}
P.~E.~T. Jorgensen and S.~Pedersen.
\newblock Orthogonal harmonic analysis of fractal measures.
\newblock {\em Electron. Res. Announc. Amer. Math. Soc.}, 4:35--42
  (electronic), 1998.

\bibitem{JP99}
P.~E.~T. Jorgensen and S.~Pedersen.
\newblock Spectral pairs in {C}artesian coordinates.
\newblock {\em J. Fourier Anal. Appl.}, 5(4):285--302, 1999.

\bibitem{KM04a}
M.~N. Kolountzakis and M.~Matolcsi.
\newblock Complex {H}adamard matrices and the spectral set conjecture.
\newblock In {\em Proc. 7th International Conference on Harmonic Analysis and
  Partial Differential Equations (El Escorial), 2004}, to appear.
\newblock (CA/11512).

\bibitem{KM04}
M.~N. Kolountzakis and M.~Matolcsi.
\newblock Tiles with no spectra.
\newblock {\em Forum Math.}, to appear.
\newblock (CA/0406127).

\bibitem{LW97}
J.~C. Lagarias and Y.~Wang.
\newblock Spectral sets and factorizations of finite abelian groups.
\newblock {\em J. Funct. Anal.}, 145(1):73--98, 1997.

\bibitem{LWWW2002}
R.~S. Laugesen, N.~Weaver, G.~L. Weiss, and E.N. Wilson.
\newblock A characterization of the higher dimensional groups associated with
  continuous wavelets.
\newblock {\em J. Geom. Anal.}, 12(1):89--102, 2002.

\bibitem{O02}
G.~{\'O}lafsson.
\newblock Continuous action of {L}ie groups on {$\R\sp n$} and frames.
\newblock {\em Int. J. Wavelets Multiresolut. Inf. Process.}, 3(2):211--232,
  2005.

\bibitem{olafsson_speegle:03}
G.~{\'O}lafsson and D.~Speegle.
\newblock Wavelets, wavelet sets, and linear actions on {$\R\sp n$}.
\newblock In {\em Wavelets, frames and operator theory}, volume 345 of {\em
  Contemp. Math.}, pages 253--281. Amer. Math. Soc., Providence, RI, 2004.

\bibitem{S}
D.~Speegle.
\newblock On the existence of wavelets for non-expansive dilation matrices.
\newblock {\em Collect. Math.}, 54(2):163--179, 2003.

\bibitem{tao:04}
T.~Tao.
\newblock Fuglede's conjecture is false in 5 and higher dimensions.
\newblock {\em Math. Res. Lett.}, 11(2-3):251--258, 2004.

\bibitem{wang:02}
Y.~Wang.
\newblock Wavelets, tiling, and spectral sets.
\newblock {\em Duke Math. J.}, 114(1):43--57, 2002.

\bibitem{WW01}
G.~Weiss and E.~N. Wilson.
\newblock The mathematical theory of wavelets.
\newblock In {\em Twentieth century harmonic analysis---a celebration (Il
  Ciocco, 2000)}, volume~33 of {\em NATO Sci. Ser. II Math. Phys. Chem.}, pages
  329--366. Kluwer Acad. Publ., Dordrecht, 2001.

\end{thebibliography}
\end{document}